\documentclass{article}\begin{document}
\centerline{\bf\large On a General Sextic Equation Solved by}
\centerline{\bf\large the Rogers-Ramanujan Continued fraction}
\vskip .4in

\centerline{Nikos Bagis}

\centerline{Stenimahou 5 Edessa}
\centerline{Pella 58200, Greece}
\centerline{bagkis@hotmail.com}
\vskip .2in

\[
\]


\[
\]

\centerline{\bf Abstract}

\begin{quote}
In this article we solve a general class of sextic equations. The solution follows if we consider the $j$-invariant and relate it with the polynomial equation's coefficients. The form of the solution is a relation of Rogers-Ramanujan continued fraction. The inverse technique can also be used for the evaluation of the Rogers-Ramanujan continued fraction, in which the equation is not now the depressed equation but another quite more simplified equation.      
\end{quote}

\section{Introductory Definitions}
We will solve the following equation  
$$
\frac{b^2}{20a}+bX^3+aX^6=C_1X^{5}\eqno{:(eq)} 
$$
or equivalent
\begin{equation}
\frac{b^2}{20a}+bX+aX^2=C_1X^{5/3}
\end{equation} 
using the $j$-invariant and the Rogers-Ramanujan continued fraction.
\[
\] 
For $\left|q\right|<1$, the Rogers Ramanujan continued fraction (RRCF) (see [2],[3],[4]) is defined as
\begin{equation}
R(q):=\frac{q^{1/5}}{1+}\frac{q^1}{1+}\frac{q^2}{1+}\frac{q^3}{1+}\cdots  
\end{equation}
From the Theory of Elliptic functions the $j$-invariant (see [5],[8]) is
\begin{equation}
j_r:=\left[\left(\frac{\eta(\frac{1}{2}\sqrt{-r})}{\eta\left(\sqrt{-r}\right)}\right)^{16}+16\left(\frac{\eta(\sqrt{-r})}{\eta\left(\frac{1}{2}\sqrt{-r}\right)}\right)^{8}\right]^3 , 
\end{equation}
where 
\begin{equation}
\eta(\tau):=e^{\pi i\tau/12}\prod^{\infty}_{n=1}\left(1-e^{2\pi i n\tau}\right)
\end{equation}
is the Dedekind's eta function and 
$$
\tau=\frac{1+\sqrt{-r}}{2} \textrm{ , } \tau=\sqrt{-r}\textrm{ , }r\textrm{ positive real. } 
$$ 
We have also in the $q$-notation
\begin{equation}
f(-q):=\prod^{\infty}_{n=1}(1-q^n) . 
\end{equation}
In what follows we use the following known result (see Wolfram pages for 'Rogers-Ramanujan Continued Fraction' and [17]):\\
If 
$$
R=R(e^{-2\pi\sqrt{r}}) , 
$$ 
then:
\begin{equation}
j_r=-\frac{\left(R^{20}-228 R^{15}+494 R^{10}+228 R^5+1\right)^3}{R^5 \left(R^{10}+11 R^5-1\right)^5}
\end{equation}
From ([3],[4]) we have 
\begin{equation}
\frac{1}{R^5(q)}-11-R^5(q)=\frac{f^6(-q)}{q f^6(-q^5)} 
\end{equation} 
The general hypergeometric function is defined as
$$
{}_pF_{q}\left[\left\{a_1,a_2,\ldots,a_p\right\},\left\{b_1,b_2,\ldots,b_q\right\},x\right]=\sum^{\infty}_{n=0}\frac{(a_1)_n(a_2)_n,\ldots,(a_p)_n}{(b_1)_n(b_2)_n,\ldots,(b_q)_n}\frac{x^n}{n!}
$$
where $(c)_n=c(c+1)\ldots (c+n-1)$, hence $(1)_n=n!$ .\\
The standard definition of the elliptic integral of the first kind (see [7],[8],[15]) is:
\begin{equation}
K(x)=\int^{\pi/2}_{0}\frac{dt}{\sqrt{1-x^2\sin^2(t)}} 
\end{equation}
\begin{equation}
K(x)=\frac{\pi}{2}{}_2F_{1}\left(\{1/2,1/2\};\{1\};x^2\right)=\frac{\pi}{2}{}_2F_{1}\left(1/2,1/2;1;x^2\right)
\end{equation}
In the notation of Mathematica we have
\begin{equation}
K(x)=\textrm{EllipticK}[x^2] 
\end{equation}
The elliptic singular modulus $k=k_r$ is defined to be the solution of the equation: 
\begin{equation}
\frac{K\left(\sqrt{1-k^2}\right)}{K(k)}=\sqrt{r} .
\end{equation}
In Mathematica's notation
\begin{equation}
k=k_r=k[r]=\textrm{InverseEllipticNomeQ}[e^{-\pi \sqrt{r}}]^{1/2} . 
\end{equation}
The complementary modulus is given by $k'_{r}=\sqrt{1-k_r^2}$. (For evaluations of $k_r$  see [5],[15],[16]).\\
Also we call $w_r:=\sqrt{k_rk_{25r}}$ noting that if one knows $w=w_r$ then (see [2]), knows $k_r$ and $k_{25r}$.

\section{Theorems}

\textbf{Proposition 1.} (see [2])\\
If $q=e^{-\pi\sqrt{r}}$ and $r$ real positive then we define
\begin{equation}
A=A_r:=\frac{f^6(-q^2)}{q^2f^6(-q^{10})}=R(q^2)^{-5}-11-R(q^2)^5
\end{equation}
then 
\begin{equation}
A_r=a_{4r}=\frac{(k_rk'_r)^2}{(w_rw'_r)^2}\left(\frac{w_r}{k_r}+\frac{w'_r}{k'_r}-\frac{w_rw'_r}{k_rk'_r}\right)^3
\end{equation}
\textbf{Theorem 1.}\\ Let $a$, $b$, $C_1$ be constants. One can solve the equation
\begin{equation}
\frac{b^2}{20a}+bX+aX^2=C_1X^{5/3},
\end{equation} 
finding $r>0$ such that
\begin{equation}
j_r=250C^3_1a^{-2}b^{-1} .
\end{equation} 
Then (15) have solution 
\begin{equation}
X=\frac{b}{250 a}A_r=\frac{b}{250a}\frac{f(-e^{-2\pi\sqrt{r}})^6}{e^{-2\pi\sqrt{r}}f(-e^{-10\pi\sqrt{r}})^6} . 
\end{equation}
\textbf{Proof.}\\
For to solve the equation (15) find $r$ such that 
\begin{equation}
j^{1/3}_r=\frac{5\cdot2^{1/3}C_1}{a^{2/3}b^{1/3}}  
\end{equation} 
Consider also the transformation of the constants
$$
3125m=\frac{b^2}{20a}\textrm{ , }250ml^{-1}=\frac{b^2}{250a}\left(b+\frac{b^2}{20a}\right)^{-1}
$$
and
$$ml^{-2}=\frac{b^2}{62500}\left(b+\frac{b^2}{20a}\right)^{-2} , $$
with inverse 
$$
l=\frac{b(20a+b)}{20a} \textrm{ , } m=\frac{b^2}{62500a} . 
$$
Then 
$$X=\frac{250m}{l(l-3125m)}x_1=\frac{250m}{l-3125m}x=\frac{b}{250a}x , $$ 
where $x_1$ satisfies
$$3125m+250x_1 m l^{-1}+x^2_1ml^{-2}=ml^{-5/3}j^{1/3}x^{5/3}_1$$
If we set $x_1=lx$, then it is
$$3125+250x+x^2=j^{1/3}x^{5/3}$$
or equivalently
\begin{equation}
3125+250A_r+A^2_r=j^{1/3}_rA^{5/3}_r
\end{equation}
Relation (19) is equivalent to equation (6), in view of (7).
Hence from Proposition 1 
$$
X=A_r=a_{4r}=\frac{b}{250a}\frac{f(-e^{-2\pi\sqrt{r}})^6}{e^{-2\pi\sqrt{r}}f(-e^{-10\pi\sqrt{r}})^6}=
$$
\begin{equation}
=\frac{b}{250a}\frac{(k_rk'_r)^2}{(w_rw'_r)^2}\left(\frac{w_r}{k_r}+\frac{w'_r}{k'_r}-\frac{w_rw'_r}{k_rk'_r}\right)^3=\frac{b}{250a}\left(R^{-5}(q^2)-11-R^5(q^2)\right)   
\end{equation}
and the proof is complete.
\[
\]
The $j$-invariant is connected with the singular modulus from the equation 
\begin{equation}
j_r=\frac{256(k^2_r+k'^4_{r})^3}{(k_rk'_r)^4} . 
\end{equation}
We can solve (21) and express $k_r$ in radicals to an algebraic function of $j_r$.\\   
The 5th degree modular equation which connects $k_{25r}$ and $k_r$ is (see [3]):
\begin{equation}
k_rk_{25r}+k'_rk'_{25r}+2\cdot4^{1/3} (k_rk_{25r}k'_rk'_{25r})^{1/3}=1
\end{equation}
We will evaluate the root of (1) first with parametrization and second with Rogers-Ramanujan continued fraction and the Elliptic-$K$ function.\\
For this it have been showed (see [19]) that if 
\begin{equation}
k_{25r}k_r=w^2_r=w^2 , 
\end{equation} 
setting the following parametrization of $w$:
\begin{equation}
w=\sqrt{\frac{L(18+L)}{6(64+3L)}} , 
\end{equation}
we get
\begin{equation} 
\frac{(k_{25r})^{1/2}}{w^{1/2}}=\frac{w^{1/2}}{(k_r)^{1/2}}=\frac{1}{2}\sqrt{4+\frac{2}{3}\left(\frac{L^{1/6}}{M^{1/6}}-4\frac{M^{1/6}}{L^{1/6}}\right)^2}+\frac{1}{2}\sqrt{\frac{2}{3}}\left(\frac{L^{1/6}}{M^{1/6}}-4\frac{M^{1/6}}{L^{1/6}}\right)
\end{equation} 
where $$M=\frac{18+L}{64+3L}$$
From the above relations we get also
\begin{equation}
-\frac{k_r-w}{\sqrt{k_rw}}=\frac{k_{25r}-w}{\sqrt{k_{25r}w}}=\sqrt{\frac{2}{3}}\left(\frac{L^{1/6}}{M^{1/6}}-4\frac{M^{1/6}}{L^{1/6}}\right)
\end{equation}
Hence we can consider the above equations as follows: Taking an arbitrary number $L$ we construct an $w$. Now for this $w$ we evaluate the two numbers $k_{25r}$ and $k_r$. Thus when we know the $w$, the $k_{r}$ and $k_{25r}$ are given from (24),(25),(26).\\ 
The result is: We can set a number $L$ and from this calculate the two inverse elliptic nome's. But we don't know the $r$. One can see (from the definition of $k_r$) that the $r$ can evaluated from equation 
\begin{equation}
r=\frac{K^2(\sqrt{1-k^2_r})}{K^2(k_r)}
\end{equation}
Hence we define
\begin{equation}
r=k^{(-1)}(x):=\frac{K^2(\sqrt{1-x^2})}{K^2(x)}
\end{equation}
However is very difficult to evaluate the $r$ in a closed form, such as roots of polynomials or else when a number $x$ is given. Some numerical evaluations indicate us that even if $x$ are algebraic numbers, (not trivial as with $k^{(-1)}\left(2^{-1/2}\right)=1$ or the cases $x=k_r$, $r=1,2,3,\dots$) the $r$ are not rational and may even not algebraics.
\[
\]
\centerline{\textbf{The algebraic representation of $X$}}
\[
\]
We know that (see [2]): 
\begin{equation}
X=X(L)=\frac{b}{250a}\frac{x_L^2(1-x^2_L)}{(w_Lw'_L)^2}\left(\frac{w_L}{x_L}+\frac{w'_L}{\sqrt{1-x^2_L}}-\frac{w_Lw'_L}{x_L\sqrt{1-x^2_L}}\right)^3
\end{equation} 
where $x_L=k_r$ is the singular modulus which corresponds to some $L$. 
\begin{equation}
C_1=\frac{a^{2/3}b^{1/3}}{5\cdot 2^{1/3}}j_{r_L}^{1/3} .
\end{equation}
The procedure is to select a number $L$ and from (24),(25) evaluate $w_L$, $x_L$ and 
\begin{equation}
w'_L=\sqrt{\sqrt{1-\frac{w_L^4}{x_L^2}}\sqrt{1-x_L^2}} . 
\end{equation}
The solution $X=X(L)$ of (1) is (29) and for this $L$ holds 
\begin{equation}
r_L=k^{(-1)}(x_L)
\end{equation}
\begin{equation}
X=\frac{b}{250a}\frac{f\left(-e^{-2\pi\sqrt{k^{(-1)}(x_L)}}\right)^6}{e^{-2\pi\sqrt{k^{(-1)}(x_L)}}f\left(-e^{-10\pi\sqrt{k^{(-1)}(x_L)}}\right)^6}
\end{equation}
\begin{equation}
j_{r_L}=250C^3_1a^{-2}b^{-1}
\end{equation}  
\begin{equation}
j_{r_L}=\frac{256(x^2_L+(1-x^2_L)^2)^3}{x^4_L(1-x^2_L)^2}=250\frac{C_1^3}{a^2b}
\end{equation}
Hence we get the next:
\[
\] 
\textbf{Theorem 2.}\\ One can find parametric solutions of (1) if for a given $L$ construct the $w_L$, $x_L$ and the complementary $w'_L$ (these values are given from (23),(24),(25),(31)). Also $x'_L=\sqrt{1-x^2_L}$. The $C_1$ must be 
\begin{equation}
C_1=\sqrt[3]{\frac{a^2 b j_{r_L}}{250}} 
\end{equation}
The solution $X=X_L$ is given from 
\begin{equation}
X=X(L)=\frac{b}{250a}\frac{x_L^2(1-x^2_L)}{(w_Lw'_L)^2}\left(\frac{w_L}{x_L}+\frac{w'_L}{\sqrt{1-x^2_L}}-\frac{w_Lw'_L}{x_L\sqrt{1-x^2_L}}\right)^3
\end{equation}
\[
\]
\textbf{Note.}\\ 
\textbf{i)} The above solution (37) works for parametric solutions (setting a $L$), as also for solutions which we know $r$, $k_r$ and $k_{25r}$. (For a related method on solving the quintic see Wolfram pages 'Quintic Equation') \\ 
\textbf{ii)} In [16] it have been shown that when one knows for some $r_0$ the $k_{r_0}$ and $k_{r_0/25}$ then can evaluate any $k_{25^nr_0}$ in radicals closed form for all $n$ positive integers. But in general the values $k_r$ and $k_{25r}$ can given from tables or with a simple PC (see [4],[5],[13],[15],[17]). 
\[
\]
\centerline{\textbf{The inverse functions method}}
\[
\]
From the analysis in [2], the solution $X$ of (1) can reduced with inverse functions as follows:\\ 
Consider the function
\begin{equation}
U(x)=\frac{256(x^2+(1-x^2)^2)^3}{x^4(1-x^2)^2} , 
\end{equation} 
The equation $U(x)=t$ have known solution with respect to $x$, which we will call $x=U^{(-1)}(t)$. Hence
\begin{equation}
\frac{256(k^2_r+(1-k^2_r)^2)^3}{k^4_r(1-k^2_r)^2}=250\frac{C_1^3}{a^2b}
\end{equation} 
or
$$
U(k_r)=250\frac{C_1^3}{a^2b}
$$
$$
k_{r}=U^{(-1)}\left(250\frac{C_1^3}{a^2b}\right)
$$
or
$$
r=k^{(-1)}\left(U^{(-1)}\left(250\frac{C_1^3}{a^2b}\right)\right)
$$
The function $k^{(-1)}(x)$ is that of (28).
\[
\]
\textbf{Theorem 3.}\\
The equation (1) have solution 
\begin{equation}
X=\frac{b}{250a}\left(R\left(e^{-2\pi\sqrt{ k^{(-1)}(\alpha)}}\right)^{-5}-11-R\left(e^{-2\pi\sqrt{k^{(-1)}(\alpha)}}\right)^{5}\right)
\end{equation}
where
\begin{equation}
\alpha=U^{(-1)}\left(250\frac{C_1^3}{a^2b}\right)
\end{equation}\[
\]
\textbf{Notes.}\\
1) Observe here that we don't need the value of $w$ and the class invariant $j$.\\ 
2) From [10] we have 
$$
R(e^{-x})=e^{-x/5}\frac{\vartheta_4(3ix/4,e^{-5x/2})}{\vartheta_4(ix/4,e^{-5x/2})},\forall x>0
$$  
Hence the solution can expressed also in theta functions. That is if $\alpha=k_r$, $r=1,2,3,\ldots$ then the solution of (1) reduced to that of evaluation of Rogers-Ramanujan continued fraction $R(q)$ with $q=e^{-\pi \sqrt{r}}$. In view of [10] we have
\footnotesize
$$
X=\frac{b}{250a}\left[e^{2\pi\sqrt{r}}\left(\frac{\vartheta_4\left(3i\pi\sqrt{r}/2,e^{-5\pi\sqrt{r}}\right)}{\vartheta_4\left(i\pi\sqrt{r}/2,e^{-5\pi\sqrt{r}}\right)}\right)^{-5}-11-e^{-2\pi\sqrt{r}}\left(\frac{\vartheta_4\left(3i\pi\sqrt{r}/2,e^{-5\pi\sqrt{r}}\right)}{\vartheta_4\left(i\pi\sqrt{r}/2,e^{-5\pi\sqrt{r}}\right)}\right)^5\right]
$$ 
\normalsize  
\\
\textbf{Example 1.}\\
The equation $$X^2+3X+\frac{9}{20}=\frac{26}{5\sqrt[3]{3}}X^{5/3}$$ have $$\alpha=U^{(-1)}\left(\frac{35152}{9}\right)=\frac{\sqrt{3}}{2}$$ hence a solution is
$$
X=\frac{3}{250}\left(R\left(e^{-2\pi\sqrt{r}}\right)^{-5}-11-R\left(e^{-2\pi\sqrt{r}}\right)^5\right)
$$
where
$$ r=\frac{K\left(\frac{1}{2}\right)^2}{K\left(\frac{\sqrt{3}}{2}\right)^2}
$$
For this $r$ the $X$ is a solution. 
\[
\]

Continuing one can set to 
\begin{equation}
X^5=\frac{b^2}{20a C_1}+\frac{b}{C_1}X^3+\frac{a}{C_1}X^6
\end{equation}
any value $X=X_0$ and $C_1=1$ then evaluate 
\begin{equation}
a=\frac{-5b X_0^3+5 X_0^5+\sqrt{5} \sqrt{4 b^2 X_0^6-10 b X_0^8+5 X_0^{10}}}{10 X_0^6}
\end{equation}
equation (42) holds always and we get that
$$
R\left(e^{-2\pi \sqrt{r}}\right)^{-5}-11-R\left(e^{-2\pi  \sqrt{r}}\right)^5=$$
$$
=\frac{25 \left(-5 b+5 X_0^2+\sqrt{5} \sqrt{4 b^2-10 b X_0^2+5 X_0^4}\right)}{b X_0^2}
$$
where $j_r=250a^{-2}b^{-1}$. 
\begin{equation}
j_r=\frac{25000 X_0^6}{b \left(-5 b+5 X_0^2+\sqrt{5} \sqrt{4 b^2-10 b X_0^2+5 X_0^4}\right)^2}   
\end{equation}   
The result is the following parametrized evaluation of the Rogers-Ramanujan continued fraction\\
\\
\textbf{Theorem 4.}
$$
A_r=R\left(e^{-2\pi \sqrt{r}}\right)^{-5}-11-R\left(e^{-2\pi  \sqrt{r}}\right)^5=
$$
\begin{equation}
=\frac{25 \left(-5b+5t^{5}+\sqrt{5} \sqrt{4 b^2-10 b t^{2}+5 t^{4}}\right)}{b}
\end{equation}
and 
\begin{equation}
j_r=\frac{25000 t^6}{b \left(-5 b+5 t^{2}+\sqrt{5} \sqrt{4 b^2-10 b t^{2}+5 t^{4}}\right)^2}   
\end{equation}
\\
\textbf{Corollary.}\\
If $$\sqrt[3]{A_r^2j_r}\textrm{ = rational}$$ then $A_r$ is of the form $$A_r=\frac{A+B\sqrt{D}}{C}$$ where $A,B,C,D$ rationals\\
\\ 
\textbf{Theorem 5.}\\   
If for a certain $r>0$ we know the value of  $R(e^{-\pi\sqrt{r}})$ in radicals, then we can evaluate both $k_r$ and $k_{25r}$ and the opposite.\\
\textbf{Proof.}\\
Suppose we know for a certain $r>0$ the value of $R(e^{-\pi\sqrt{r}})$, (the correspondence between  $R(e^{-\pi\sqrt{r}})$ and $R(e^{-2\pi\sqrt{r}})$ is given by (96) bellow or see [13]). Then from (6) we know the value of $j_r$ and from (21) we know $k_r$. Let also $q=e^{-\pi\sqrt{r}}$, $r>0$ and $v_r=R(q)$, then it have been proved by Ramanujan that
$$
v^5_{r/25}=v_r\frac{1-2v_r+4v_r^2-3v_r^3+v_r^4}{1+3v_r+4v_r^2+2v_r^3+v_r^4} , 
$$ 
Hence we can get the value of $R(e^{-\pi\sqrt{r}/5})$.
Hence again form (6) we find $j_{r/25}$ and from (21) the value of $k_{r/25}$. But from relation (53) bellow knowing $k_r$ and $k_{r/25}$ we can evaluate all $k_{25^nr}$, $n=1,2,\ldots$ and consequently $k_{25r}$ as a special case. 
The oposite follow from Proposition 1.\\ 
\\
\textbf{Theorem 6.}\\
The solution $U_0$ of the equation
\begin{equation}
U_0=j_r^{1/3}\left(125-\sqrt{12500+U_0}\right)^{5/3}
\end{equation}
is
$$ 
U_0=U(j_r)=\sum^{\infty}_{n=1}\frac{j^{n/3}_r}{n!}\left[\frac{d^{n-1}}{da^{n-1}}\left(125-\sqrt{12500+a}\right)^{5n/3}\right]_{a=0}
$$
If $X=x_0$ is root of 
$$
U_0=X^2+250X+3125
$$
then 
$$
3125+250x_0+x_0^2=j_r^{1/3}(-1)^{1/3}x_0^{5/3}\eqno{(47a)}
$$
and 
$$
x_0=X=X_r=A_r=R(e^{-2\pi\sqrt{r}})^{-5}-11-R(e^{-2 \pi \sqrt{r}})^5\eqno{(47b)}
$$
\textbf{Proof.}\\
Consider (1), then make the change of variable $U_0=X^2+250X+3125$, we arrive to (47). The Legendre inversion theorem states that the solution of $y=af(y)$ (see [7] pg.132-133) is 
$$
y=\sum^{\infty}_{n=1}\frac{a^n}{n!}\left[\frac{d^{n-1}}{dx^{n-1}}f(x)^n\right]_{x=0}
$$
\\
This theorem works for $j_r$ small, for example with $j_1=1728$ it converges very slowly but for $r$ such that $j_r=800$, ($r$-complex) we get numerical evaluations and hence also theoretical.\\
\\
\textbf{Theorem 7.}\\
If
$$
c_n:=\left[\frac{d^{n-1}}{da^{n-1}}\left(125-\sqrt{12500+a}\right)^{5n/3}\right]_{a=0}
$$
then 
$$
c_n=\frac{5^6}{3}(-1)^{n+1}n\cdot10^{-5n/3}{}_2F_{1}\left[\frac{5n}{6},\frac{5n+3}{6};\frac{2(n+3)}{3};\frac{1}{5}\right]\frac{\Gamma(5n/3)}{\Gamma(2+2n/3)} 
$$
\textbf{Proof.}\\
Recall a theorem of Euler (see [18] pg.306-307). If the root of 
$$
aqx^p+x^q=1
$$ 
is $x$, then 
$$
x^n=\frac{n}{q}\sum^{\infty}_{k=0}\frac{\Gamma(\{n+pk\}/q)(-q a)^k}{\Gamma(\{n+pk\}/q-k+1)k!}
$$  
Hence from the fact that 
$$
w=\frac{250(125-\sqrt{12500+x})}{3125-x}
$$ 
is solution of 
$$
acb^{-2}w^2+w=1\textrm{ , }a=1\textrm{ , }b=-250\textrm{ , }c=-3125+x 
$$
we get
\begin{equation}
\left(-125+\sqrt{12500+x}\right)^n=\frac{n}{250^n}\sum^{\infty}_{k=0}\frac{\Gamma(n+2k)(-1)^k}{\Gamma(n+k+1)62500^k k!}(-3125+x)^{k+n}
\end{equation}
Using the formula
$$
f^{(\nu)}(x_0)=\sum^{\infty}_{n=0}\frac{f^{(\nu+n)}(x_0)}{n!}(-x_0)^n
$$
the result follows.\\
\\ 
Theorem 7 is for numerical evaluations since the  hypergeometric series can more easily computed than the $(n-1)$th derivative of the $5n/3$ power of $125-\sqrt{12500+x}$.\\
\\
\textbf{Corollary.}\\
For every 'suitable' value of $x_0$ such that $X_r=x_0 :(a)$, $X_r$ is of (47b) exists a $r$ solution of (a) such that
$$
X_r^2+250X_r+3125=
$$
$$
=3^{-1}\cdot5^6\sum^{\infty}_{n=1}(-1)^{n+1}n\frac{ \Gamma(5n/3)}{\Gamma(2+2n/3)}{}_2F_{1}\left[\frac{5n}{6},\frac{5n+3}{6};\frac{2(n+3)}{3};\frac{1}{5}\right] \frac{(10^{-5}j_r)^{n/3}}{n!}
$$  
\\
\textbf{Example 2.}\\
For $X_r=x_0=-12$, we have $r=-0.186710441...-i0.251574161...$
and 
$$
X_r^2+250X_r+3125=269=\sum^{\infty}_{n=1}\frac{(-j_r)^{n/3}}{n!}\left[\frac{d^{n-1}}{dz^{n-1}}\left(125-\sqrt{12500+z}\right)^{5n/3}\right]_{z=0}
$$
\\  
\textbf{Example 3.}\\
Consider the equation 
$$
X^2+250X+3125=2(-1)^{1/3}10^{2/3}X^{5/3} 
$$
Then clearly $j=j_r=800$ and a solution is
$$
X=x_0=-125+\sqrt{12500+\sum^{\infty}_{n=1}\frac{(2 \sqrt[3]{100})^{n}}{n!}\left[\frac{d^{n-1}}{dz^{n-1}}\left(125-\sqrt{12500+z}\right)^{5n/3}\right]_{z=0}}
$$ 

\section{Applications}

\textbf{Example 4.}\\
Set $L=1/3$ then $$w_N=w_1(L)=w_1\left(\frac{1}{3}\right)=\frac{1}{3}\sqrt{\frac{11}{78}}$$ 
and 
$$k_N=x_L=x_1(L)=x_1\left(\frac{1}{3}\right)=$$
$$=\frac{\frac{1}{3}\sqrt{\frac{11}{78}}}{\left(\frac{-4(\frac{11}{13})^{1/6}+(\frac{13}{11})^{1/6}}{\sqrt{6}}+\frac{1}{2}\sqrt{4+\frac{2}{3}\left(-4\left(\frac{11}{13}\right)^{1/6}+\left(\frac{13}{11}\right)^{1/6}\right)^2}\right)^2}$$
and
$$k_{25N}=x_2(L)=x_2\left(\frac{1}{3}\right)=$$
$$=\frac{1}{3}\sqrt{\frac{11}{78}}\left(\frac{-4(\frac{11}{13})^{1/6}+(\frac{13}{11})^{1/6}}{\sqrt{6}}+\frac{1}{2}\sqrt{4+\frac{2}{3}\left(-4\left(\frac{11}{13}\right)^{1/6}+\left(\frac{13}{11}\right)^{1/6}\right)^2}\right)^2$$
where the $N$ is given by
$$
N=r_L=r_{1/3}=\frac{K^2\left(\sqrt{1-x_1\left(\frac{1}{3}\right)^2}\right)}{K^2\left(x_1\left(\frac{1}{3}\right)\right)}
$$
From the value of $x_L$ we obtain $j_{r_L}$ and hence the corresponding $C_1$ in radicals-closed form and hence $X=X_L$ from (37) and (31). The numbers $a$, $b$ take arbitrary values.\\
We note here that in future application of this method one must tabulate values of $(r,w_r)$ and not $j_r$ or $k_r$ which follow from these of $w_r$. This can be done in some cases using the Main Theorem in [16] and the solution (37) of Theorem 2 of the present paper.     
\[
\]
Form [16] we have if
\begin{equation}
Q(x)=\frac{\left(-1-e^{\frac{1}{5}y}+e^{\frac{2}{5} y}\right)^5}{\left(e^{\frac{1}{5}y}-e^{\frac{2}{5}y}+2 e^{\frac{3}{5} y}-3 e^{\frac{4}{5}y}+5 e^{y}+3 e^{\frac{6}{5}y}+2 e^{\frac{7}{5}y}+e^{\frac{8}{5} y}+e^{\frac{9}{5}y}\right)}
\end{equation}
$$
y=\textrm{arcsinh}\left(\frac{11+x}{2}\right)
$$
\begin{equation}
Y=U_0(X)=\sqrt{-\frac{5}{3 X^2}+\frac{25}{3 X^2 h(X)}+\frac{X^4}{h(X)}+\frac{h(X)}{3 X^2}}
\end{equation}
$$
h(x)=\left(-125-9x^6+3\sqrt{3}\sqrt{-125x^6-22x^{12}-x^{18}}\right)^{1/3}
$$
\begin{equation}
U_1(Y)=X=\sqrt{-\frac{1}{2 Y^2}+\frac{Y^4}{2}+\frac{\sqrt{1+18 Y^6+Y^{12}}}{2 Y^2}} .
\end{equation}
and
\begin{equation}
P(x)=P[x]=U_0[Q^{1/6}[U_1[x]]]\textrm{ and }
P^{(n)}(x)=(P\underbrace{\circ\ldots\circ}_{n} P)(x)  
\end{equation}
then
\begin{equation}
k_{25^nr_0}=\sqrt{1/2-1/2\sqrt{1-4\left(k_{r_0}k'_{r_0}\right)^2\prod^{n}_{j=1}P^{(j)}\left[\sqrt[12]{\frac{k_{r_0}k'_{r_0}}{k_{r_0/25}k'_{r_0/25}}}\right]^{24}}}
\end{equation}
\textbf{Example 5.}
$$
k_{1/5}=\sqrt{\frac{9+4 \sqrt{5}+2 \sqrt{38+17 \sqrt{5}}}{18+8 \sqrt{5}}}
$$
$$
k_{5}=\sqrt{\frac{9+4 \sqrt{5}-2 \sqrt{38+17 \sqrt{5}}}{18+8 \sqrt{5}}}
$$
\begin{equation}
k_{125}=\sqrt{\frac{1}{2}-\frac{1}{2}\sqrt{1-(9-4\sqrt{5})P[1]^2}}  
\end{equation}
\textbf{Example 6.}\\
It is
$$k_1=\frac{1}{\sqrt{2}}$$
$$
k_{25}=\frac{1}{\sqrt{2 \left(51841+23184 \sqrt{5}+12 \sqrt{37325880+16692641 \sqrt{5}}\right)}}
$$
Hence
$$
k_{625}k'_{625}=\frac{1}{2 \left(161+72 \sqrt{5}\right)}P\left[161-72\sqrt{5}\right]
$$
and hence
\begin{equation}
k_{625}=\sqrt{\frac{1}{2}-\frac{1}{2}\sqrt{1-\left(\frac{P\left[161-72\sqrt{5}\right]}{161+72 \sqrt{5}}\right)^2}}
\end{equation}
By this way we can evaluate every $k_{r}$ which is $r=4^l9^m25^nr_0$ when $k_{r_0}$ and $k_{r_0/25}$ are known, $l,m,n\in\bf N\rm$. 
\[
\]
\textbf{Note.}
In the case that  
\begin{equation}
\left(\frac{L}{M}\right)^{1/6}=A=\frac{f}{g}=3^{a}\frac{2p+1}{2h+1} , 
\end{equation}
where $f$, $g$ positive integers, with $f<g$ and $p, h\neq 0(mod4)$, ($a\in\bf Z\rm-\left\{-1,0,1\right\}$), we can find $w$ from a given $x=k_r$ which is of the form 
$$
x=k_r=\frac{t_1\sqrt{t_2}}{\left(t_3+\sqrt{t_4}\right)^2}
$$
where $t_i$, $i=1,2,3,4$ rationals.
\\
In view of (25) and the action of the command recognize, (which is needed to put number $x$ into his form) the output will be an octic equation with step 2 containing nested square roots:  
$$
<<NumberTheory'Recognize'
$$
$$
\textrm{Solve[Reduce[N[x,1000],16,v]==0,v]}  
$$
The smallest root it will be 
$$\sqrt{D}=\sqrt{4096+88\left(3^{a}\frac{2p+1}{2h+1}\right)^6+\left(3^{a}\frac{2p+1}{2h+1}\right)^{12}}$$
One can see that for these $x$'s the $f$ and $g$ are given from the Diofantie equation
\begin{equation}
9D=g^{12}\left(4096+88\frac{f^6}{g^6}+\frac{f^{12}}{g^{12}}\right)
\end{equation}
hence the number $A$ will be known and
\begin{equation}
w^2=\frac{4096-20 A^6+A^{12}+\left(-64+A^6\right) \sqrt{4096+88 A^6+A^{12}}}{108 A^6}
\end{equation}
Hence we have the value of $X$ in radicals.
\[
\]
If for example
$$    
x=k_r=\frac{-37754085 \sqrt{3}+3 \sqrt{476791023769787}}{77 \sqrt{2} \left(-435+\sqrt{224799}\right)^2}
$$
then for all $C_1$, $a$, $b$ such that
$$
j_{r}=\frac{256(x^2+(1-x^2)^2)^3}{x^4(1-x^2)^2}=250\frac{C_1^3}{a^2b}
$$
with Mathematica and the package 'Recognize' we evaluate
$$
Solve[Recognize[N[x,1000],16,v]==0,v]
$$
which gives the value of $x$ in the desired form.
The solution that corresponds to $x$ have smallest square root 
$$\sqrt{D}=\sqrt{1430373071309361}$$
the command 'Reduce' give us the $f$ and $g$
$$
Reduce[9D==(4096+88(f/g)^6+(f/g)^{12})g^{12},\left\{f,g\right\},Integers]
$$
Hence we get the values $f=7$, $g=11$ and $w$. The solution (29) is
$$
X=\frac{A}{35153041^3}[5579801448-11724990 \sqrt{224799}+
$$
$$
+\sqrt{6362897839 \left(9487950991-20011160 \sqrt{224799}\right)}]^3
$$
where
$$
A=A_1+B_1-\frac{1}{2}\sqrt{A_2+B_2\sqrt{224799}}
$$
\footnotesize
$$
A_1=\frac{93573266991461291403517623659291588}{29148873138738228269392700625}
$$
$$
B_1=\frac{572443990137 \sqrt{15689769027087558724590525868357405909071}}{22336301255738105953557625}
$$
$$
A_2=\frac{70047382179949155201445598248892391809778079434882871423635278599557376}{849656805258255011387371218620407164266412150030875390625}
$$
$$
B_2=\frac{12578872085638673940246389099496155002610886190948544424950380544}{72342001299127714890367919848480814326642158367890625}
$$ 
\normalsize
\\
\textbf{Theorem 8.}\\ 
When $(L/M)^{1/6}$ is rational, we can always find $k_{25r}$ from $k_r$. 
\\
\\
\textbf{Example 7.}\\
Set $a=(k_rk'_r)^2$, $b=(ww')$, then equation (1) have solution $$X=\frac{m_5^3}{250} , $$
where $m_5$ is the multiplier (see [3]).\\ 
Hence 
$$
3125\frac{(ww')^4}{(k_rk'_r)^2}+250(ww')^2m^3_5+(k_rk'_r)^2m_5^6=(k_rk'_r)^{4/3}(ww')^{2/3}j^{1/3}_rm^5_5
$$
\[
\]
\textbf{Example 8.}\\
We will find a solution of the equation
\begin{equation}
\frac{3125}{16}+125X+4X^2=132X^{5/3}
\end{equation}
in radicals using Theorem 3.\\
\textbf{Solution.}\\
It is $a=4$, $b=125$, and we have to solve $j_r=287496$ or equivalently $r=4$.\\
Hence a solution of (59) is:
\begin{equation}
X=\frac{125e^{4\pi}}{250\cdot4}\frac{f(-e^{-4\pi})^6}{f(-e^{-20\pi})^6}=\frac{1}{8}\left(R(e^{-4\pi})^{-5}-11-R(e^{-4\pi})^5\right)
\end{equation} 
The exact root in radicals can be found but is very large and complicated with our method. We give a way how one can obtain it:\\   
It is known that
\begin{equation}
R(e^{-2\pi})=-\frac{1+\sqrt{5}}{2}+\sqrt{\frac{5+\sqrt{5}}{2}}
\end{equation}
But from the duplication formula (see [4],[13]):\\
If $u=R(q)$ and $\nu=R(q^2)$, then
\begin{equation}
\frac{\nu-u^2}{\nu+u^2}=u\nu^2 . 
\end{equation}
Hence we find the value of $R(e^{-4\pi})$ in radicals and hence the solution of (59) using (60),(61),(62).\\
The root using the program Mathematica is  
$$
X=\frac{143375}{16}+\frac{64125 \sqrt{5}}{16}+\frac{1}{2} \sqrt{\frac{20553203125}{32}+\frac{9191671875 \sqrt{5}}{32}}
$$
In this case it is more convenient to use Mathematica's command Solve. But in other cases these solutions can not found.
\[
\]
From the above result we have shown that
$$
\frac{1}{8}\left(R(e^{-4\pi})^{-5}-11-R(e^{-4\pi})^5\right)=
$$
$$
=\frac{143375}{16}+\frac{64125 \sqrt{5}}{16}+\frac{1}{2} \sqrt{\frac{20553203125}{32}+\frac{9191671875 \sqrt{5}}{32}} .
$$
One can see that if we set $$Y_{\tau}:=\frac{b}{250a}\left(R\left(e^{2\pi i \tau}\right)^{-5}-11-R\left(e^{2  \pi i \tau}\right)^5\right)\eqno{:(a)}$$
Then if $$\tau=\frac{1+\sqrt{-r}}{2} \textrm{ or } \tau=\sqrt{-r} , $$ 
where $r$ positive integer in some cases we can evaluate $Y_{\tau}$ solving directly the equation (1), with parameters $a=4$, $b=125$ and $C_1$ depended on $j_{\tau}$.\\ Some examples are 
\footnotesize
$$
\frac{1}{16}\left[R\left(e^{\pi(i-\sqrt{51})}\right)^{-5}-11-R\left(e^{\pi(i-\sqrt{51})}\right)^5\right]=
$$
$$
=-\frac{125}{16} [5541103+1343914 \sqrt{17}+\sqrt{61407604829690+14893531819350 \sqrt{17}}+
$$
$$
+4 \surd \{\frac{1}{6140760482969+1489353181935 \sqrt{17}}(94272348104055803848937570+$$
$$+22864402871059934148609270 \sqrt{17}+$$
$$+2063164169063100077 \sqrt{170 \left(6140760482969+1489353181935 \sqrt{17}\right)}+$$
\begin{equation}
+8506643792036854023 \sqrt{61407604829690+14893531819350 \sqrt{17}})\}] . 
\end{equation}
\[
\]
\normalsize
\begin{equation}
Y_{\sqrt{-1/5}}=\frac{5\sqrt{5}}{8} . 
\end{equation}
\begin{equation}
Y_{\sqrt{-2/5}}=\frac{5}{8}\left(5+2\sqrt{5}\right) . 
\end{equation}
\begin{equation}
Y_{\sqrt{-3/5}}=\frac{5}{16} \left(25+11 \sqrt{5}\right) . 
\end{equation}
\begin{equation}
Y_{\sqrt{-4/5}}=\frac{5}{16} \left(25+13 \sqrt{5}+5 \sqrt{58+26 \sqrt{5}}\right) . 
\end{equation}
\begin{equation}
Y_{\sqrt{-5/5}}=\frac{125}{8} \left(2+\sqrt{5}\right) . 
\end{equation}
\begin{equation}
Y_{\sqrt{-6/5}}=\frac{5}{8} \left(50+35 \sqrt{2}+3 \sqrt{5 \left(99+70 \sqrt{2}\right)}\right) . 
\end{equation}
\begin{equation}
Y_{\sqrt{-9/5}}=\frac{5}{8} \left(225+104 \sqrt{5}+10 \sqrt{1047+468 \sqrt{5}}\right) . 
\end{equation}
\begin{equation}
Y_{\sqrt{-12/5}}=\frac{5}{16} \left(1690+975 \sqrt{3}+29 \sqrt{6755+3900 \sqrt{3}}\right) . 
\end{equation}
\begin{equation}
Y_{\sqrt{-14/5}}=\frac{5}{8} \left(1850+585 \sqrt{10}+7 \sqrt{5 \left(27379+8658 \sqrt{10}\right)}\right).
\end{equation}
\begin{equation}
Y_{\sqrt{-17/5}}=\frac{5}{8}\left(5360+585 \sqrt{85}+4 \sqrt{3613670+391950 \sqrt{85}}\right).
\end{equation}
We describe the method bellow.\\
For some $r$ positive rational we find the value of $j_{r/5}$; this can be done with the command 'Recognize' of the program Mathematica (if $j_{r/5}$ is root of a small degree algebraic polynomial equation). Then we find $C_1$ (from (16)) and for the values $a=4$, $b=125$ there will be $$Y_{\tau}=\textrm{root of equation (1)} . $$
In many cases of such $r$, equation (1) can solved in radicals with Mathematica (we have not find the reason yet), but still in others not. Hence we get relations like (63)$-$(73).

\section{More Theorems and Results}

\textbf{Theorem 9.} (Conjecture)\\ For every positive real $r$, we have
\begin{equation}
Y_{\sqrt{-r/5}}Y_{\sqrt{-r^{-1}/5}}=\frac{125}{64} .
\end{equation}
\[
\]
If $l$, $m$, $t$ and $d$ are integers and  
\begin{equation}
Y_{\sqrt{-r/5}}=\frac{l+m\sqrt{d}}{t}
\end{equation}
then
\begin{equation}
l^2-m^2d=t^2\frac{125}{64}
\end{equation}
In general we conjecture that
\[
\]
\textbf{Theorem 10.} (Conjecture)\\ 
If $r=a_1/b_1$ with $a_1, b_1\in\bf N\rm$ and $\textrm{GCD}(a_1,5)=1$, $\textrm{GCD}(b_1,5)=1$ then 
\begin{equation}
deg\left(Y_{\sqrt{-r/5}}\right)=deg\left(j_{\sqrt{-r/5}}\right)
\end{equation} 
\[
\]
For example if $deg\left(Y_{\sqrt{-r/5}}\right)=4$, then
\begin{equation}
Y_{\sqrt{-r/5}}=A+B\sqrt{D}
\end{equation}
where $deg(A)=deg(D)=2$ and 
\begin{equation}
A^2-B^2 D=\frac{125}{64}U
\end{equation}
where $deg(U)=2$ or $U=1$. If $U\neq1$ then $U=l+m\sqrt{d}$ and also if $j_{\sqrt{-r/5}}$ have smallest nested square root $\sqrt{d}$, then $UU^{*}=l^2-m^2d=1$. The symbol $*$ denotes the algebraic conjugate.\\ 
Hence for example if $r=6$ then $d=2$ and
$$
j_{\sqrt{-6/5}}=8640 [25551735275-18067805280 \sqrt{2}-$$
$$-196 \sqrt{10 \left(3399058140008707-2403497060447490 \sqrt{2}\right)}]
$$
then $U=l+m\sqrt{2}$ with 
$$
l^2-2m^2=1 . 
$$ 
We solve the above Pell's equation. The solution we looking for, taking the smallest to higher order solutions, for this example with $r=6$ is $l_1=99$ and $m_1=70$. Hence  $A^2-B^2D=\frac{125}{64}(99+70\sqrt{2})$.\\
Now we assume that $A=k_1+l_1\sqrt{d}$, again with $d=2$ and $D=k_2+l_2\sqrt{d}$, etc... We proceed solving Pell's equations.\\
\\
\textbf{Theorem 11.}\\ 
For a given $r\in \bf N\rm$ and $deg\left(Y_{\sqrt{-r/5}}\right)=2$, $4$, or $8$, if the smallest nested root of $j_{\sqrt{-r/5}}$ is $\sqrt{d}$ then we can evaluate the Rogers-Ramanujan continued fraction with integer parameters.\\ 
\textbf{i)} In the case $deg\left(Y_{\sqrt{-r/5}}\right)=2$ then
\begin{equation}
Y_{\sqrt{-r/5}}=\frac{l+m\sqrt{d}}{t}
\end{equation}
where 
\begin{equation}
l^2-m^2d=1\textrm{ and } l,m,d\in\bf N\rm
\end{equation}
\textbf{ii)} In the case $deg\left(Y_{\sqrt{-r/5}}\right)=4$ we have\\
\textbf{a)} If $U\neq\frac{125}{64}$, then
\begin{equation}
Y_{\sqrt{-r/5}}=\frac{5}{8}\sqrt{\left(a_0+\sqrt{-1+a_0^2}\right)}\left(\sqrt{5+p}-\sqrt{p}\right) 
\end{equation}
where 
\begin{equation}
Y_{\sqrt{-r/5}}Y^{*}_{\sqrt{-r/5}}=\frac{125}{64}\left(a_0+\sqrt{a_0^2-1}\right) , 
\end{equation}
with, $a_0$ positive integer, is solution of $l^2-m^2d=1$. Hence $l=a_0$ and $m=d^{-1/2}\sqrt{a_0^2-1}$ is positive integer. The parameter $p$ is positive rational can be found from the numerical value of $Y_{\sqrt{-r/5}}$.\\
\textbf{b)} If $U=\frac{125}{64}$, then
\begin{equation}
Y_{\sqrt{-r/5}}=A+\frac{1}{8}\sqrt{-125+64A^2} ,
\end{equation}
where we set $A=k+l\sqrt{d}$. Then a starting point for the evaluation of the integers $k$, $l$ will be the relation 
\begin{equation}
l^2=\frac{(A-k)^2}{d}=\textrm{ square of integer }
\end{equation}
\textbf{iii)} If $deg\left(Y_{\sqrt{-r4^{-1}5^{-1}}}\right)=4$, then we can evaluate $Y_{\sqrt{-r5^{-1}}}$.\\
It holds $deg\left(Y_{\sqrt{-r5^{-1}}}\right)=8$, the minimal polynomial of $Y_{\sqrt{-r5^{-1}}}/Y_{\sqrt{-r4^{-1}5^{-1}}}$ is of degree 4 or 8 and symmetric. Hence it can be reduced in at most 4th degree polynomial, hence it is solvable. Thus it remains the evaluation of $Y_{\sqrt{-r4^{-1}5^{-1}}}$, which can be done with the help of step (ii). 
\begin{equation}
Y_{\sqrt{-r5^{-1}}}=\frac{5}{8}\sqrt{a_0+\sqrt{-1+a^2_0}}\left(\sqrt{p+5}-\sqrt{p}\right)2^{-1}\left(\sqrt{x+4}-\sqrt{x}\right)  
\end{equation}
where $x=a_1+b_1\sqrt{d}+c\sqrt{a_2+b_2\sqrt{d}}$, $a_1$, $b_1$, $a_2$, $b_2$, $c$ integers and
$$ Y_{\sqrt{-r5^{-1}4^{-1}}}=\frac{5}{8}\sqrt{a_0+\sqrt{-1+a^2_0}}\left(\sqrt{p+5}-\sqrt{p}\right)
$$
\[
\] 
\textbf{Example 9.}\\
For $r=68=4\cdot17$ and from (73) we have $d=85$
$$
x=a_1+b_1\sqrt{85}+c\sqrt{a_2+b_2\sqrt{85}}
$$ 
$$
Y_{\sqrt{-68/5}}/Y_{\sqrt{-17/5}}=2^{-1}\left(\sqrt{x+4}-\sqrt{x}\right)
$$
$$
a_1=2891581250 \textrm{, } b_1=313636050 \textrm{, } c=12960
$$
$$ 
a_2=99557521554 \textrm{, } b_2=10798529365  
$$
hence 
$$
Y_{\sqrt{-68/5}}=Y_{\sqrt{-17/5}}2^{-1}\left(\sqrt{x+4}-\sqrt{x}\right)=
$$
$$
=\frac{5}{16}\left(5360+585 \sqrt{85}+4\sqrt{3613670+391950 \sqrt{85}}\right)\left(\sqrt{x+4}-\sqrt{x}\right)
$$
\\
\textbf{Theorem 12.}\\
If $r=a_1/b_1$ with $deg(j_{r/5})=\nu\leq 4$, then $deg(A_{r/5})=\nu$ and equation (1) (with $a$, $b$ rationals) can solved in radicals.\\
\\
\textbf{Application.}\\
If $r=3/4$ then $deg(j_{3/20})=4$ and $A_{r/5}$ is solution of
$$ 
15625-2112500 v+443375 v^2-16900 v^3+v^4=0   
$$
hence
$$
A_{3/20}=R\left(e^{-\pi\sqrt{3/5}}\right)^{-5}-11-R\left(e^{-\pi\sqrt{3/5}}\right)^5=
$$
$$
=\frac{5}{2} \left(1690-975 \sqrt{3}+29 \sqrt{6755-3900 \sqrt{3}}\right)
$$
\\
\textbf{Theorem 13.}\\ 
If $Q(x):=x^5$ then
\footnotesize
\begin{equation}
\frac{1}{j^{1/3}_{\tau}}\left[R\left(e^{2\pi i \tau}\right)^{-5}-11-R\left(e^{2\pi i\tau}\right)^5\right]^{1/3}=
\sqrt[3]{\frac{-125}{j_{\tau}}+\sqrt{\frac{12500}{j^2_{\tau}}{+Q\left(\sqrt[3]{\frac{-125}{j_{\tau}}+\ldots}\right)}}}
\end{equation}
\normalsize
\textbf{Proof.}\\
Equation (1) for $a=1$, $b=250j^{-1}_{\tau}$, $C_1=1$ can be written in the form 
\begin{equation}
(X^3-a_1)^2-b_1=X^5+c_1 , 
\end{equation}
where $a_1=-125j^{-1}_{\tau}$, $b_1=12500j^{-2}_{\tau}$, $c_1=0$\\
Hence $Y_{\tau}$ we can be expressed in nested periodical functions. This completes the proof.  
\[
\]
\textbf{Example 10.}\\
If 
$$
C^3_1=32a^2b
$$
then
$$
X=\frac{b}{250a}\left(R(e^{-2\pi\sqrt{2}})^{-5}-11-R(e^{-2\pi\sqrt{2}})^{5}\right)
$$
\[
\]
\centerline{\textbf{Equation (1) and the Derivative of Rogers-Ramanujan Continued fraction}}
\[
\]
From [11] it is known that if 
\begin{equation}
N(q)=q^{5/6}f(-q)^{-4}\frac{R'(q)}{R(q)}
\end{equation}
and $N(q^2)=u(q)=u$, $N(q^3)=h(q)=h$ and $N(q)=v(q)=v$, then 
\begin{equation}
5u^6-u^2v^2-125u^4v^4+5v^6\stackrel{?}{=0}
\end{equation}
and
\begin{equation}
125h^{12}+h^3v^3+1125h^9v^3+1125h^3v^9+1953125h^9v^9-125v^{12}\stackrel{?}{=0}
\end{equation}
which are solvable. But from [12] we have
\begin{equation}
\frac{5R'(q)}{R(q)\left(R(q)^{-5}-11-R(q)^5\right)^{1/6}}=f^4(-q)q^{-5/6}
\end{equation}
or 
\begin{equation}
N(q)=\frac{1}{5}\left(R(q)^{-5}-11-R(q)^5\right)^{1/6}
\end{equation} 
Hence the solution of (1) can also given in the form  
\begin{equation}
X=X_r=\frac{125b}{2a}N(q^2)^6
\end{equation}
and from (85) we have 
\begin{equation}
2^{2/3} a^{1/3} b^{1/3} (X_rX_{4r})^{1/3}+\frac{10\cdot 2^{1/3} a}{b^{1/3}} (X_r X_{4r})^{2/3}-2 a^{2/3} (X_r+X_{4r})=0
\end{equation}
\textbf{Note.}\\
One can prove relation (90) using (89),(93) and the duplication formula (see [13]):
\begin{equation}
\frac{R(q^2)-R^2(q)}{R(q^2)+R^2(q)}=R(q)R^2(q^2)
\end{equation}
The same method can work and with other higher modular equations of the derivative but the evaluations are very difficult even for a program.
\[
\]
Another interesting note that can simplify the problem is the singular moduli of the fifth base (see [14],[15]):
\begin{equation}
u(x)={}_2F_1\left(\frac{1}{6},\frac{5}{6};1;x\right) . 
\end{equation}
In this case we have
\begin{equation}
j_r=\frac{432}{\beta_r(1-\beta_r)}=\frac{250C_1^3}{a^2b} , 
\end{equation}
where $\beta_r$ is the solution of
\begin{equation}
\frac{{}_2F_1\left(\frac{1}{6},\frac{5}{6};1;1-\beta_r\right)}{{}_2F_1\left(\frac{1}{6},\frac{5}{6};1;\beta_r\right)}=\sqrt{r}
\end{equation}
The moduli $\beta_r$ can evaluated from $k_r$ and the opposite from the relation
\begin{equation}
\frac{256(k^2_r+(1-k^2_r)^2)^3}{k^4_r(1-k^2_r)^2}=\frac{432}{\beta_r(1-\beta_r)}
\end{equation}
\[
\]
\textbf{Proposition 2.}\\
The equation (1) have solution 
\begin{equation}
X=\frac{b}{250a}\left(R\left(e^{-2\pi\sqrt{ \beta^{(-1)}(\alpha)}}\right)^{-5}-11-R\left(e^{-2\pi\sqrt{\beta^{(-1)}(\alpha)}}\right)^{5}\right) , 
\end{equation}
where 
$$\alpha=\frac{1}{2}-\sqrt{\frac{1}{4}-\frac{216 a^2b}{125 C_1^3}} \textrm{ , }  \beta^{(-1)}(x)=\left(\frac{u(1-x)}{u(x)}\right)^2 . 
$$
\\
\textbf{Corollary.}\\
The equation 
$$
aX^2+bX+\frac{b^2}{20a}=\frac{6a^{2/3}b^{1/3}}{5\beta_r^{1/3}(1-\beta_r)^{1/3}}X^{5/3} 
$$
admits solution $X=A_r$.
\newpage

\centerline{\bf References}\vskip .2in

\noindent

[1]: M.Abramowitz and I.A.Stegun. 'Handbook of Mathematical Functions'. Dover Publications, New York. 1972.

[2] Nikos Bagis. 'Parametric Evaluations of the Rogers-Ramanujan Continued Fraction'. International Journal of Mathematics and Mathematical Sciences. Vol (2011) 

[3]: B.C. Berndt. 'Ramanujan`s Notebooks Part III'. Springer Verlag, New York (1991)

[4]: B.C. Berndt. 'Ramanujan's Notebooks Part V'. Springer Verlag, New York, Inc. (1998)

[5]: D. Broadhurst. 'Solutions by radicals at Singular Values $k_N$ from New Class Invariants for $N\equiv3\;\; mod\;\; 8$'. arXiv:0807.2976 (math-phy).        

[6]: I.S. Gradshteyn and I.M. Ryzhik. 'Table of Integrals, Series and Products'. Academic Press (1980).

[7]: E.T.Whittaker and G.N.Watson. 'A course on Modern Analysis'. Cambridge U.P. (1927)

[8]: J.V. Armitage W.F. Eberlein. 'Elliptic Functions'. Cambridge University Press. (2006)

[9]: Bruce C. Berndt and Aa Ja Yee. 'Ramanujans Contributions to Eisenstein Series, Especially in his Lost Notebook'. (page stored in the Web). 

[10]: Nikos Bagis and M.L. Glasser. 'Jacobian Elliptic Functions, Continued Fractions and Ramanujan Quantities'. arXiv:1001.2660v1 [math.GM] 2010.

[11]: Nikos Bagis. 'Generalizations of Ramanujans Continued fractions'. arXiv:11072393v2 [math.GM] 7 Aug 2012.  

[12]: Nikos Bagis, M.L. Glasser. 'Integrals related with Rogers-Ramanujan continued fraction and $q$-products'. arXiv:0904.1641. 10 Apr 2009.

[13]: Bruce C. Berndt, Heng Huat Chan, Sen-Shan Huang, Soon-Yi Kang, Jaebum Sohn and Seung Hwan Son. 'The Rogers-Ramanujan Continued Fraction'. (page stored in the Web).

[14]: N.D. Bagis and M.L. Glasser. 'Conjectures on the evaluation of alternative modular bases and formulas  approximating $1/\pi$'. Journal of Number Theory. (Elsevier), (2012).

[15]: J.M. Borwein and P.B. Borwein. 'Pi and the AGM'. John Wiley and Sons, Inc. New York, Chichester, Brisbane, Toronto, Singapore. (1987)

[16]: Nikos Bagis. 'Evaluation of Fifth Degree Elliptic Singular Moduli'. arXiv:1202.6246v1. (2012)

[17]: W. Duke. 'Continued Fractions and Modular Functions'. (page stored in the Web)

[18]: B.C. Berndt. 'Ramanujan`s Notebooks Part I'. Springer Verlag, New York. (1985)

[19]: Nikos Bagis. 'The complete evaluation of Rogers Ramanujan and other continued fractions with elliptic functions'. arXiv:1008.1340v1. (2010)

\end{document}